\documentclass[a4paper,english,frenchb,11pt]{amsart}
\usepackage[latin1]{inputenc}
\usepackage{amssymb,babel}\NoAutoSpaceBeforeFDP
\usepackage[T1]{fontenc}
\usepackage{aeguill,calrsfs}

\usepackage{ifpdf}
\ifpdf\usepackage[pdftex]{hyperref}\else\usepackage[hypertex]{hyperref}\fi

\usepackage{amscd}

\newtheorem{theo}{Théorème}
\newtheorem{coro}[theo]{Corollaire}
\newtheorem{lemm}[theo]{Lemme}
\theoremstyle{remark}
\newtheorem{rema}[theo]{Remarque}

\renewcommand{\leq}{\leqslant}
\renewcommand{\geq}{\geqslant}

\newcommand{\setC}{\mathbb{C}}
\newcommand{\setR}{\mathbb{R}}
\newcommand{\cC}{\mathcal{C}}
\newcommand{\cI}{\mathcal{I}}
\newcommand{\cF}{\mathcal{F}}
\newcommand{\cL}{\mathcal{L}}
\newcommand{\cN}{\mathcal{N}}
\newcommand{\cO}{\mathcal{O}}
\newcommand{\cT}{\mathcal{T}}
\newcommand{\db}{\bar{\partial}}

\DeclareMathOperator{\Res}{Res}

\begin{document}

%\date{\today}
\author{Olivier Biquard}
\address{IRMA, Université Louis Pasteur et CNRS, 7 rue René Descartes,
  F-67084 STRASBOURG CEDEX}
\email{olivier.biquard@math.u-strasbg.fr}
\title{Sur les variétés CR de dimension 3 et les twisteurs}
\selectlanguage{english}
\begin{abstract}
  We prove that any real analytic strictly pseudoconvex CR 3-manifold
  is the boundary (at infinity) of a unique selfdual Einstein metric
  defined in a neighborhood. The proof uses a new construction of
  twistor space based on singular rational curves.
\end{abstract}
\selectlanguage{frenchb}
\maketitle

\section*{Introduction}
Dans l'article fondateur \cite{LeB82}, LeBrun montre qu'une variété
conforme $X^3$, analytique réelle, est l'infini conforme d'une
métrique d'Einstein autoduale, unique dans un voisinage de $X$.
Ici, infini conforme signifie que la métrique d'Einstein admet un pôle
d'ordre deux au bord, mais sa structure conforme s'étend au bord $X$
pour induire la structure conforme initiale.

Une question similaire peut être posée en partant d'une variété CR
strictement pseudoconvexe $X^3$. Le but de cet article est d'y répondre
positivement. L'exemple standard de la théorie est l'espace
hyperbolique complexe, dont un modèle est la métrique de Bergmann sur
la boule de $\setR^4$,
\begin{equation}
  \label{eq:Bergmann}
  g_B = \frac{1}{1-r^2}\sum_1^4 dx_i^2
        + \frac{r^2}{(1-r^2)^2}(dr^2+r^2\eta^2)
\end{equation}
où $\eta$ est la 1-forme sur la sphère $S^3$ obtenue comme forme de
connexion du fibré de Hopf $S^3\to S^2$.

Plus généralement, si $X^3$ est munie d'une distribution de contact $H$,
noyau de la forme de contact $\ker \eta$, et d'une structure CR $J$ sur
$H$, de sorte que
\begin{equation}
  \label{eq:gamma}
  \gamma(x,y)=d\eta(x,Jy)
\end{equation}
soit une métrique sur $H$, on dira que $X$ est l'infini conforme d'une
variété riemannienne $(M^4,g)$ si $X=\partial M$, et, en choisissant une
fonction $x$ sur $M$ s'annulant exactement sur $X$ à l'ordre un, on a
le comportement asymptotique suivant de $g$ près du bord $X$:
\begin{equation}
  \label{eq:ach}
  g \sim \frac{dx^2+\eta^2}{x^2}+\frac{\gamma}{x} .
\end{equation}
Cette condition impose que le tenseur de courbure de $g$ soit
asymptotique à celui de la métrique de Bergmann, voir \cite{Biq00}.
Un changement conforme $(\eta,\gamma)\to(f\eta,f\gamma)$ par une fonction $f>0$ ne modifie
pas ce comportement asymptotique, puisqu'on retrouve (\ref{eq:ach}) en
faisant $x\to\frac{x}{f}$.

Les exemples explicites de métriques d'Einstein autoduales avec une
variété CR comme infini conforme sont les métriques $SU_2$-invariantes
sur la boule \cite{Hit95}, et des métriques toriques \cite{CalSin04}.

D'un autre côté existe une construction plus abstraite : une structure
CR sur la sphère $S^3$, proche de la structure standard, est l'infini
conforme d'une unique métrique d'Einstein sur la boule, complète,
petite déformation de la métrique de Bergmann \cite{Biq00}. Parmi ces
structures CR, celles dont le remplissage d'Einstein est autodual
forment une famille de dimension infinie, dont l'espace tangent en la
structure standard est le noyau de l'opérateur différentiel $(\db
\sharp \db)^*$ agissant sur l'espace des sections du fibré $\Omega^{0,1}\otimes
T^{1,0}$, espace vu comme l'espace tangent à toutes les structures CR
(cette condition, démontrée dans \cite{Biq05}, est reformulée ainsi
dans l'article de synthèse \cite{Biq06}).  Comme cette condition infinitésimale
est de nature locale, cela conduit naturellement à la question
suivante : le problème, pour une structure CR, d'être l'infini
conforme d'une métrique d'Einstein autoduale, admet-il une obstruction
de nature locale ? ou, au contraire, la condition infinitésimale
indiquée plus haut est-elle bien une obstruction globale au remplissage par
une métrique \emph{complète} sur la boule ?

Dans le cas réel, le théorème de LeBrun évoquée plus haut montre qu'il
n'y a pas d'obstruction locale, et la condition pour remplir une
métrique conforme sur la sphère par une métrique d'Einstein autoduale
complète sur la boule, trouvée dans \cite{Biq02}, est en effet de
nature globale (phénomène de «~fréquences positives~» prédit par
LeBrun)---même la condition infinitésimale est globale.

Le résultat suivant indique qu'il en est de même dans le cas
complexe, en montrant l'absence d'obstruction locale.

\begin{theo}\label{th:1}
  Soit $X$ une variété CR strictement pseudoconvexe de dimension 3,
  analytique réelle. Alors $X$ est l'infini conforme d'une métrique
  d'Einstein autoduale $g$ définie dans un petit voisinage de $X$ ($g$
  est analytique réelle et a le comportement (\ref{eq:ach}) près de
  $X$).  Le germe de $g$ le long de $X$ est unique.
\end{theo}

En particulier, le développement asymptotique d'une métrique
d'Einstein autoduale près de son infini conforme est complètement
déterminé, et la métrique est analytique réelle. Cela contraste
fortement avec le problème similaire pour une métrique
Kähler-Einstein d'infini conforme une structure CR : dans le
développement formel apparaissent des termes indéterminés et des termes
logarithmiques \cite{Fef76,LeeMel82}. 

Remarquons que la formule (\ref{eq:Bergmann}) pour la métrique de
Bergmann garde un sens pour $r>1$, et définit une métrique
Kähler-Einstein autoduale, de signature (2,2), sur le fibré en disques
holomorphes $\cO(1)$ au-dessus de la droite projective $P^1$. La
métrique $g$ construite dans le théorème précédent est analytique
réelle, et en la prolongeant de l'autre côté de $X$, on obtient le
résultat suivant.

\begin{coro}
  Sous les mêmes hypothèses, la variété $X$ est l'infini conforme
  d'une métrique d'Einstein autoduale de signature (2,2), unique,
  définie dans un petit voisinage de $X$, de comportement
  asymptotique donné par la formule (\ref{eq:ach}) avec $x<0$.
\end{coro}

La transformation twistorielle de Penrose (voir \cite{AtiHitSin78})
est un outil puissant dans la construction de métriques autoduales
d'Einstein.  Ces deux résultats sont démontrés à l'aide d'une nouvelle
construction twistorielle, où le bord est réalisé comme ensemble de
courbes rationnelles ayant un point double.

Décrivons le plan de l'article.  On commence par décrire dans la
section \ref{sec:lesp-des-twist} l'espace des twisteurs de l'espace
hyperbolique complexe, motivant la construction twistorielle
proprement dite dans les sections \ref{sec:les-twisteurs} et
\ref{sec:fch}. Celle-ci consiste à relever, dans la complexification
de $X$, les feuilletages en courbes holomorphes induits par les fibrés
$T^{1,0}$ et $T^{0,1}$ de la structure CR, à un fibré de formes
différentielles, puis à recoller les deux espaces de feuilles.  La
construction rappelle celle des twisteurs du cotangent d'une variété
kählérienne par Feix \cite{Fei01}, dans laquelle sont recollés des
espaces de fonctions affines sur les sous-variétés intégrales des
distributions $T^{1,0}$ et $T^{0,1}$ dans la complexification de la
variété. On réalise alors la variété initiale comme espace de courbes
rationnelles dans l'espace des twisteurs, mais ici s'ajoute une difficulté
technique : ces courbes sont singulières.  Après lissification
(section \ref{sec:cou-rat}), on accomplit dans les sections
\ref{sec:vcs} et \ref{sec:la-transf-twist} la transformée twistorielle
inverse pour déterminer le comportement de la métrique ainsi
construite.

\noindent\emph{Remerciements.} Je remercie Olivier Debarre pour
m'avoir longuement expliqué le problème de lissification des courbes
rationnelles singulières.

\section{L'espace des twisteurs de l'espace hyperbolique complexe}\label{sec:lesp-des-twist}

On considère l'espace $\setC^{1,2}$. Dans des coordonnées complexes
$(z^1,z^2,z^3)$, on a donc une forme hermitienne
$|z^1|^2-|z^2|^2-|z^3|^2$. On notera $\cF$ l'espace des drapeaux de
$\setC^3$, consistant d'une droite $D$ et d'un plan $P\supset D$. L'espace
$\cF$ dispose de deux projections sur $P^2_\setC$, données
respectivement par $(D,P)\to D$ et $(D,P)\to P$. En chaque point, la somme
des espaces tangents à chacune des fibres des projections est un
sous-fibré holomorphe du fibré tangent de codimension un---en fait une
structure de contact holomorphe.

L'espace hyperbolique complexe $H^2_\setC$ se décrit comme l'espace
des droites positives de $\setC^{1,2}$. Son espace des twisteurs,
$\cN$, est le domaine de $\cF$ décrit par
$$ \cN = \{(D,P)\in\cF, D<0, P \text{ de signature }(1,1) \}. $$
La projection twistorielle $\pi:\cN\to H^2_\setC$ est
$$ \pi(D,P)=D^\perp \cap P , $$
sa structure réelle
$$ \tau(D,P)=(P^\perp,D^\perp), $$
sa structure de contact holomorphe est la restriction à $\cN$ de celle
de $\cF$, elle est en chaque point transverse à la fibre de $\pi$.

Le bord de $H^2_\setC$, la sphère $S^3$, s'identifie à l'espace des
droites isotropes de $\setC^{1,2}$. La projection twistorielle s'étend
au-dessus du bord, et l'espace des twisteurs $\cT=\cN|_{S^3}$
au-dessus de $S^3$ a deux composantes $\cT_±$, et s'écrit
\begin{equation}
  \cT = \cT_+ \cup_{S^3} \cT_- .\label{eq:T-mod}
  \end{equation}
En effet, ces deux composantes sont décrites par:
\begin{itemize}
\item $\cT_+$ est l'espace des $(D,P)$ tels que $D$ est isotrope et
  $P\supset D$ ; sur cette composante $\pi(D,P)=D$;
\item $\cT_-$ est l'espace des $(D,P)$ tels que $P$ soit isotrope et
  $D\subset P$ ; sur cette composante $\pi(D,P)=P^\perp$.
\end{itemize}
Les deux composantes s'intersectent en les couples $(D,P)$, où $D$ et
$P$ sont isotropes, avec $D=P^\perp=\pi(D,P)$ ; cette intersection est
donc égale à $S^3$. La structure réelle $\tau$ échange $\cT_+$ et
$\cT_-$ en fixant $S^3$.

Enfin, explicitons complètement la famille de courbes rationnelles de
$\cF$ dans laquelle se trouvent les fibres de la projection twistorielle
$\pi$. Cette famille est de dimension 4, et est paramétrée par un couple
$(d,p)$ d'une droite $d$ et d'un plan $p$ de $\setC^3$ ; la courbe
correspondante dans $\cF$ est
\begin{equation}
  C(d,p) = \{(D,P), D\subset p, P=D+d \} .\label{eq:Cdp}
\end{equation}
Cette courbe est lisse si $d\not\subset p$ ; si au contraire $d\subset p$,
c'est-à-dire $(d,p)\in\cF$, alors les courbes dégénèrent vers 
$$ C(d,p) = \{(D,p), D\subset p \} \cup \{(d,P), P\supset d\} , $$
qui est une courbe rationnelle avec un point double $(d,p)$.

\subsubsection*{Généralisation}
Ce comportement du bord de l'espace des twisteurs de l'espace
hyperbolique complexe est en réalité commun à toutes les métriques
avec le comportement asymptotique (\ref{eq:ach}). En effet, l'espace
des twisteurs d'une variété conforme $X$ de dimension 3 s'identifie au
projectivisé des directions nulles du cotangent $(\Omega^1_X)^\setC$,
et la forme de contact holomorphe est héritée de la forme de Liouville
de $(\Omega^1_X)^\setC$, voir \cite{LeB84}. Appliquons cette
description dans le cas suivant : si $X$ est une variété CR
strictement pseudoconvexe, avec repère local
$(\eta,\theta^1,\theta^{\bar 1})$ de $(\Omega^1_X)^\setC$, où $\eta$
est une forme de contact, $\theta^1$ engendre l'espace $T^{1,0}\subset
(\ker \eta)\otimes\setC$ de la structure CR, et
$d\eta=i\theta^1\land\theta^{\bar 1}$, alors l'espace des twisteurs
conformes pour la métrique $\frac{dx^2+\eta^2}{x^2}+\frac{\gamma}{x}$
s'identifie au-dessus d'un point $p\in X$ à la courbe d'équation
\begin{equation}
 x a^2 + bc = 0\label{eq:courbe-singuliere-x}
\end{equation}
dans $P(\Omega^1_pX)^\setC$, où une 1-forme est paramétrée par $a\eta+b\theta^1+c\theta^{\bar 1}$.
Quand $x$ tend vers 0, la courbe tend vers la courbe singulière
d'équation
\begin{equation}
 bc = 0 ,\label{eq:courbe-singuliere}
 \end{equation}
donc le bord de l'espace des twisteurs peut être considéré comme
constitué de courbes rationnelles singulières. Dans la suite, nous
construirons un espace des twisteurs coïncidant au bord avec la
réunion de ces courbes singulières.

\section{Les twisteurs}\label{sec:les-twisteurs}

Partons d'une variété $X^3$ satisfaisant les hypothèses du théorème
\ref{th:1}. La structure CR est déterminée par le fibré en droites
complexes $T^{1,0}\subset T^\setC_X$ des vecteurs de type $(1,0)$, tel qu'en
posant $T^{0,1}=\overline{T^{1,0}}$, on ait $T^{1,0}\cap T^{0,1}=0$.

La variété $X$ admet un fibré holomorphe tangent de rang 2,
$T'=T^\setC_X/T^{0,1}$, muni de l'opérateur $\db$ défini par
$$ \db_X Y = [X,Y] $$
pour $X\in T^{0,1}$ et $Y\in T'$. Cet opérateur agit aussi sur le dual
$\Omega'=(T^{0,1})^\perp \subset (T^\setC_X)^*$. En particulier, si $\eta$ est
une 1-forme réelle s'annulant sur $T^{1,0}\oplus T^{0,1}$, alors
$\db \eta\in \Omega^{0,1}\otimes\Omega'$ s'identifie à la forme de Levi $d\eta$. Puisque
$X$ est strictement pseudoconvexe, en tout point $-i d\eta(\xi,\bar{\xi})>0$
pour $\xi\in T^{1,0}$ non nul.

Au moins localement, la variété $X$ admet une complexification
$X^\setC$, munie d'une structure réelle $\tau$. Les distributions
$T^{1,0}$ et $T^{0,1}$ s'y étendent comme sous-fibrés holomorphes de
$T_{X^\setC}$, et y sont échangées par la structure réelle. Nous noterons
$C^+_x$ et $C^-_x$ les courbes intégrales de $T^{0,1}$ et
$T^{1,0}$ passant par $x$. Les fibrés $T'$ et $\Omega'$ se prolongent aussi
à $X^\setC$.

Dans la suite, nous nous restreindrons toujours à un petit ouvert $U$
de $X^\setC$ autour d'un point $x_0\in X$, stable par $\tau$. Les
constructions, canoniques, se recollent dans un petit voisinage de $X$.

L'opérateur $\db$ du fibré $\Omega'$ sur $X$ se prolonge sur $X^\setC$ en une
dérivation holomorphe, $\nabla$, du fibré $\Omega'$ le long des courbes
$C^+$.  Passant au projectivisé $P\Omega'$, les directions horizontales de
$\nabla$ définissent un feuilletage de dimension 1 de $P\Omega'$, dont on notera
$\cN^+$ l'espace des feuilles, de dimension 3.

Symétriquement, on dispose d'une dérivation le long des courbes
$C^-$ sur le fibré $\Omega''=(T^{1,0})^\perp=\overline{\Omega'}$. Les directions
horizontales du projectivisé $P\Omega''$ définissent à nouveau un
feuilletage de dimension 1 de $P\Omega''$, dont on notera $\cN^-$ l'espace
des feuilles. Il est clair que $\cN^+$ et $\cN^-$ sont échangés par
la structure réelle $\tau$. Quitte à restreindre l'ouvert $U$, on peut
supposer $\cN^±$ lisses.

Les multiples de la forme de contact $\eta$ vivent aussi bien dans $\Omega'$
que dans $\Omega''$. On obtient ainsi des applications
$i_±:X^\setC\to\cN^±$ en associant à un point $x$ la classe $[\eta_x]$ de
la forme $\eta_x$ dans $\cN^±$.
\begin{lemm}
  Les applications $i_±:X^\setC\to\cN^±$ sont des immersions.
\end{lemm}
\begin{proof}
  Vérifions le pour $i_+$. Il est clair que la différentielle de $i_+$
  est injective dans les directions transverses aux courbes $C_+$.
  D'un autre côté, le long des courbes $C_+$, la dérivée $\nabla\eta$ d'une
  forme de contact $\eta$ s'identifie à $d\eta$, donc, en tout point, est
  en dehors de $\Omega^{0,1}\otimes\setC \eta$. La section de $P\Omega'$ définie par
  $\eta$ n'est ainsi parallèle en aucun point, donc la différentielle de
  $i_+$ le long de $C_+$ est injective.
\end{proof}

À nouveau en restreignant l'ouvert $U$ de $X^\setC$, on peut supposer
que les $i_±$ sont des plongements de $U$ comme ouvert de $\cN^±$.
\begin{lemm}
  On peut choisir le petit ouvert $U$ de sorte que l'application
  $U\to\cN^+×\cN^-$ définie par $x\to(i_+(x),i_-(x))$ soit propre.
\end{lemm}
\begin{proof}
  Près du point $x_0\in X$ fixé, quitte à choisir $U$ assez petit,
  l'espace $\cF^+$ des courbes $C^+$ et l'espace $\cF^-$ des courbes
  $C^-$ sont lisses, et l'application $\phi:U\to\cF^+×\cF^-$, définie par
$$ \phi(x) = (C_x^+,C_x^-) , $$
  est une immersion. Quitte à restreindre encore $U$, l'image de $\phi$
  dans $\cF^+×\cF^-$ est une sous-variété. Fixons à présent un petit
  ouvert $V$ de $\cF^+×\cF^-$ contenant $\phi(x_0)$, et posons
  $U=\phi^{-1}(V)$. Il est maintenant clair que
$$ \phi : U \longrightarrow V $$
  est propre.
  
  On dispose de deux projections $p_±:\cN^±\to\cF^±$, et
  $\phi(x)=(p_+i_+(x),p_-i_-(x))$. Le lemme résulte alors immédiatement
  de la propreté de $\phi$.
\end{proof}

Le lemme nous indique que le recollement de $\cN^+$ et $\cN^-$ en
identifiant les ouverts $i_+(U)$ et $i_-(U)$,
\begin{equation}
  \label{eq:def-N}
  \cN = \cN^+ \amalg_U \cN^- ,
\end{equation}
est une variété. La structure réelle $\tau$ agit sur $\cN$ en échangeant
$\cN^+$ et $\cN^-$, avec points fixes $U\cap X$.

\section{La forme de contact holomorphe}
\label{sec:fch}

Au-dessus de l'espace total du fibré $P\Omega'$ sur $X$ on dispose d'un
fibré en droites complexes $L=\cO(1)$. On peut identifier la fibre de
$L$ en un point $\alpha\in \Omega'_x$ à l'espace
\begin{equation}
  T'_x/\alpha^\perp .\label{eq:L-alpha}
\end{equation}
Le fibré $L$ descend en un fibré, que nous noterons encore $L$, sur
l'espace des feuilles $\cN_+$. De même sur $P\Omega''$ nous disposons d'un
fibré $\cO(1)$ défini à partir de $T''$ par la même formule
(\ref{eq:L-alpha}). Au point $(x,[\eta])$ de $\cN_+$ ou $\cN_-$, les
fibrés $L_x$ s'identifient tous deux à $TX^\setC/(T^{1,0}\oplus T^{0,1})$
et se recollent donc pour donner un fibré $L$ sur $\cN$.

L'espace total du fibré $\Omega' \subset T^*_{X^\setC}$ hérite de la 1-forme de
Liouville $\lambda$ de $T^*_{X^\setC}$. Le long d'une courbe $C^+$, son noyau
contient par définition le vecteur horizontal $X_H$ associé à $\nabla$, et
en réalité
\begin{equation}
 i_{X_H}d\lambda=0 . \label{eq:iXdl}
\end{equation}
(Cette équation pourrait servir de définition du vecteur horizontal
$X_H$.) Comme $\lambda$ s'annule sur les fibres du fibré $\Omega' \to X^\setC$,
elle descend sur le projectivisé $P\Omega'$ comme 1-forme à valeurs dans
$L$, puis par (\ref{eq:iXdl}) sur l'espace des feuilles $\cN^+$,
toujours à valeurs dans $L$. On a donc démontré la première partie du
lemme suivant.

\begin{lemm}
  La distribution $\ker \lambda$ sur $P\Omega'$ descend en une structure de
  contact holomorphe sur $\cN^+$. De même on obtient une structure de
  contact holomorphe sur $\cN^-$.
\end{lemm}
\begin{proof}
  Il reste à montrer que la distribution est de contact. Notons $p$ la
  projection $P\Omega' \to X^\setC$. Si $\xi\in P\Omega'_x$, l'espace tangent
  $T_{[\xi]}\cN^+=T_\xi P\Omega'/\setC X_H$ peut se représenter comme
  $p^{-1}(T')=T'_x\oplus T_\xi(p^{-1}(x))$, où on a choisi un relèvement de
  $T'$ dans $T_\xi\Omega'$. Dans cette écriture, la forme $d\lambda$ est le
  crochet de dualité, elle est non nulle sur $T'_x\otimes T_\xi(p^{-1}(x)) \subset T'_x\otimes\Omega'_x$.
\end{proof}

\begin{lemm}
  Le plongement $i_±:X^\setC\to\cN^±$ ramène la structure de contact de
  $\cN^±$ sur la structure de contact standard $T^{1,0}\oplus T^{0,1}$ de $X^\setC$.
\end{lemm}
\begin{proof}
  C'est évident à partir de la formule $i_±(x)=[\eta_x]$.
\end{proof}

\begin{coro}
  Les structures de contact de $\cN^+$ et $\cN^-$ se recollent
  au-dessus de $X^\setC$ pour définir une forme de contact $\eta^c$ sur
  $\cN$, à valeurs dans le fibré en droites $L$.\qed
\end{coro}

\section{Les courbes rationnelles}
\label{sec:cou-rat}

Pour chaque point $x\in X^\setC$, les images dans $\cN$ des
projectivisés $P\Omega'_x$ et $P\Omega''_x$ définissent deux courbes
rationnelles $\cC^+_x$ et $\cC^-_x$ de $\cN$, qui se coupent au point
$[\eta_x]$ et y sont transverses. On définit la courbe
$$ \cC_x = \cC^+_x \cup \cC^-_x . $$
Les $(\cC_x)_{x\in U}$ forment une famille à trois paramètres de $P^1$
avec points doubles ; l'involution réelle $\tau$ transforme $\cC_x$ en
$\cC_{\tau(x)}$ en échangeant les parties positives et négatives.
Au-dessus des points réels, c'est-à-dire au-dessus de $X$, la famille
des $\cC_x$ dessine un espace de twisteurs qui généralise le bord de celui de
l'espace hyperbolique complexe, décrit par (\ref{eq:T-mod}), et
donne le résultat attendu (\ref{eq:courbe-singuliere}) pour une
métrique avec comportement asymptotique similaire.

\begin{lemm}
  Chaque courbe $\cC^±_x$ a fibré normal $\cO\oplus\cO$.
\end{lemm}
\begin{proof}
  On a le diagramme
$$\begin{CD}
  P\Omega' @> \pi > > \cN^+ \\
  @VV p V \\
  X^\setC
\end{CD}$$
Soit $x\in X^\setC$, le fibré normal de $\cC^+_x=\pi(P\Omega'_x)$ dans $\cN^+$
s'identifie au fibré $T_{P\Omega'}/\pi^*T_{\cC^+_x}$ sur
$P\Omega'_x$. Or $\pi^*T_{\cC^+_x}=p^*(T^{0,1})$, d'où sur $P\Omega'_x$ l'égalité
$$T_{P\Omega'}/\pi^*T_{\cC^+_x}=p^*(T')=\cO\oplus\cO . $$
\end{proof}

\begin{lemm}\label{lem:defC11}
  Les déformations d'une courbe $\cC_x$ ($x\in U$) forment une famille
  $(\cC_y)_{y\in V}$ paramétrée par une variété $V$ de dimension 4, et
  contenant $U$. Dès que $y\in V-U$, la courbe rationnelle $\cC_y$ est
  lisse, de fibré normal $\cO(1)\oplus\cO(1)$.
\end{lemm}

Pour démontrer le lemme, nous commençons par un résultat général, qui
m'a été expliqué par O.~Debarre.

\begin{lemm}\label{lem:defC}
  Soit $\cC$ une courbe rationnelle avec un point double dans une
  variété complexe $\cN^3$, dont chaque composante a fibré normal
  $\cO\oplus\cO$. Alors l'espace des déformations de $\cC$ est une
  famille lisse de dimension 4, dans laquelle :
  \begin{itemize}
  \item les courbes à point double forment une sous-famille lisse de
    dimension 3 ;
  \item les courbes lisses ont fibré normal $\cO(1)\oplus\cO(1)$ ou $\cO(2)\oplus\cO$.
  \end{itemize}
\end{lemm}

Les deux fibrés normaux possibles indiqués par le lemme se réalisent
effectivement. Le cas $\cO(1)\oplus\cO(1)$ est celui des twisteurs de
l'espace hyperbolique complexe décrit dans la section
\ref{sec:lesp-des-twist}. Le cas $\cO(2)\oplus\cO$ se réalise dans
l'exemple suivant : la courbe $\cC_0=(P^1×\{y\}×\{z\})\cup
(\{x\}×P^1×\{z\})$ dans $\cN=P^1×P^1×P^1$ satisfait les hypothèses du
lemme \ref{lem:defC}, la famille des déformations contient des courbes
lisses $\cC$ incluses dans la surface $S=P^1×P^1×\{z\}\subset\cN$, dont le fibré
normal est $N_{\cC/S}\oplus N_{S/ \cN}=\cO(2)\oplus\cO$.

\begin{proof}[Démonstration du lemme \ref{lem:defC}]
  On part d'une courbe singulière $\cC=\cC_1\cup \cC_2$, avec point
  singulier $p$. Puisque le fibré normal de $\cC_1$ est trivial, les
  déformations de $\cC_1$ forment une famille de dimension 2 de courbes
  feuilletant un voisinage ouvert $W$ de $\cC_1$ : en particulier, par
  un point de $W$ passe exactement une courbe de la famille.  La même
  chose étant valable pour $\cC_2$, les déformations de $\cC$ conservant
  un point double sont paramétrées par le point double, libre de
  bouger dans un voisinage de $p$. Donc les déformations avec point
  double forment une famille lisse de dimension 3.

  Passons aux déformations générales : on notera $f$ (resp. $f_i$)
  l'injection de $\cC$ (resp. $\cC_i$) dans $\cN$. On calcule la
  cohomologie du faisceau $f^*T_\cN$ sur $\cC$ : sur chaque composante
  $\cC_i$, on a 
  \begin{equation}
    f_i^*T_\cN=\cO\oplus\cO\oplus\cO(2),\label{eq:fiTN}
      \end{equation}
  et la suite exacte
$$ 0 \to f_1^* T_\cN(-p) \to f^* T_\cN \to f_2^* T_\cN \to 0 $$
  conduit immédiatement à 
  \begin{align*}
    H^0(\cC,f^* T_\cN)&=H^0(\cC_1,f_1^*T_\cN(-1))\oplus H^0(\cC_2,f_2^*T_\cN)=\setC^7\\ H^1(\cC,f^* T_\cN)&=0.
  \end{align*}
De l'annulation du $H^1$ on déduit (voir \cite{Kol96} ou
\cite[proposition 4.24]{Deb01}) que $f$ admet des désingularisations
$g:P^1\to\cN$. En outre, la famille des morphismes $g$ est lisse, de
dimension 7. Divisant par le groupe des paramétrisations de $P^1$, de
dimension 3, il reste une famille de déformations de $\cC$ de
dimension 4.

  Reste à trouver le fibré normal d'une courbe lisse $g:P^1\to\cN$ de
  la famille. Puisque $H^1(\cC,f^* T_\cN)=0$, on a par semi-continuité
  $H^1(P^1,g^*T_\cN)=0$. Calculons de même $H^1(P^1,g^*T_\cN(-1))$.
  Sur la courbe singulière $\cC$, prenons un faisceau $\cL$ de degré
  $0$ sur $C_1$ et $-1$ sur $C_2$, alors la suite exacte
$$ 0 \to f_1^* T_\cN(-p) \to f^* T_\cN\otimes\cL \to f_2^* T_\cN(-1) \to 0 $$
  conduit, comme ci-dessus, à l'annulation
$$ H^1(\cC,f^*T_\cN\otimes\cL) = 0 , $$
  qui implique, par semi-continuité, $H^1(P^1,g^*T_\cN(-1))=0$. Les
  possibilités pour le fibré $g^*T_\cN$ sont alors minces : écrivons
$$ g^*T_\cN = \cO(a) \oplus \cO(b) \oplus \cO(c) , \quad a \leq b \leq c , $$
  les entiers $a$, $b$ et $c$ satisfont les conditions suivantes :
  \begin{enumerate}
  \item $a+b+c=4$ car le degré de $g^*T_\cN$ est égal au degré de
    $f^*T_\cN$ sur $\cC$, donc à la somme des degrés de $f_i^*T_\cN$
    sur chaque composante $\cC_i$ ; vu (\ref{eq:fiTN}), chacun de ces
    degrés est égal à 2 ;
  \item $a\geq 0$ car $H^1(P^1,g^*T_\cN(-1))=0$ ;
  \item $g^*T_\cN$ admet le sous-faisceau $T_{P^1}=\cO(2)$, avec quotient
    localement libre, ce qui impose $c=2$, ou bien $b=2$ et $c>2$.
  \end{enumerate}
Il n'y a que deux solutions à ces contraintes :
$$ a=0, b=c=2, \qquad \text{ou bien} \qquad a=b=1, c=2. $$
Le lemme est démontré.
\end{proof}

\begin{proof}[Démonstration du lemme \ref{lem:defC11}]
  Le lemme \ref{lem:defC} fournit une famille de courbes rationnelles
  déformant la courbe singulière $\cC_x$, mais il reste à comprendre
  le fibré normal. Pour cela, nous allons, en utilisant les
  dilatations du groupe de Heisenberg, approcher la structure CR sur
  $X^3$ par la structure standard de la sphère $S^3$.

  Rappelons brièvement le modèle du groupe de Heisenberg $H$. Dans des
  coordonnées $(\sigma,u)\in \setR×\setC$, la forme de contact s'écrit
  \begin{equation}
    \eta = d\sigma + \Im(\bar{u}du) ,\label{eq:f-ctct}
      \end{equation}
  et la structure CR est définie par l'espace $T_0^{1,0}$ engendré par
  le vecteur $\frac{\partial}{\partial u}+\frac{i\bar{u}}{2}\frac{\partial}{\partial\sigma}$. Toute
  la structure est invariante sous les dilatations paraboliques 
  \begin{equation}
    \phi_t(\sigma,u)=(t^2\sigma,tu) .\label{eq:dil-H}
  \end{equation}
  Bien entendu, cette structure n'est autre que celle de la sphère
  standard $S^3$, privée d'un point. D'ailleurs les dilatations $\phi_t$
  s'étendent en des isométries de l'espace hyperbolique complexe :
  dans le groupe résoluble $\setR_+ \ltimes H$,
  \begin{equation}
    \phi_t(s,\sigma,u)=(ts,t^2\sigma,tu) ;\label{eq:dil-R}
  \end{equation}
  la coordonnée $s\geq0$ peut être considérée comme une coordonnée
  s'annulant exactement à l'ordre 1 sur le bord. Nous noterons $\cN_0$
  l'espace des twisteurs correspondant, décrit section
  \ref{sec:lesp-des-twist} et aussi plus loin section \ref{sec:vcs}.

  Revenons à présent à la structure CR sur $X$. Par le lemme de
  Darboux, on peut trouver autour du point $x\in X$ des coordonnées
  $(\sigma,u)\in \setR×\setC$ autour de l'origine, telles que la forme de
  contact soit donnée par la formule (\ref{eq:f-ctct}), et de plus la
  structure CR de $X$, déterminée par son espace $T^{1,0}$,
  coïncide à l'origine avec celle du groupe de Heisenberg. Alors, les
  dilatations paraboliques $\phi_t$ ramènent, quand $t\to0$, la structure
  CR de $X$ vers celle du groupe de Heisenberg :
  $$ \phi_t^*T^{1,0} \longrightarrow T_0^{1,0} . $$
  
  Nous noterons $\cN_t=\phi_t^*\cN$ l'espace des twisteurs pour la
  structure CR $\phi_t^*T^{1,0}$, au-dessus d'un petit ouvert fixé de
  $X$ autour de $x$ : via $\phi_t$, les espace $\cN_t$ s'identifient à
  des voisinages de plus en plus petits de la courbe singulière
  $\cC_x$ dans $\cN$, et convergent vers $\cN_0$ quand $t$ tend vers $0$.
  
  Fixons une petite boule, $V$, dans la famille de courbes
  rationnelles déformant la courbe singulière $\cC_x$ de $\cN$, et
  $U\subset V$ la sous-famille de dimension 3 des courbes singulières (qui
  s'identifie à un petit voisinage de $x$ dans $X^\setC$). Nous
  conviendrons que l'origine dans $V$ correspond à la courbe $\cC_x$.
  Les dilatations $\phi_t$ agissent sur $V$ en préservant $U$ sur lequel
  elles agissent par les dilatations paraboliques (\ref{eq:dil-H}). En
  fait, on peut paramétrer l'espace $V$ de courbes rationnelles de
  $\cN_t$ de sorte que l'action de $\phi_t$ sur $V$ soit donnée par la
  formule (\ref{eq:dil-R}) dans des coordonnées \emph{complexes}
  $(s,\sigma,u,\bar{u})$ sur $V$, où $u$ et $\bar{u}$ sont considérées
  comme coordonnées complexes indépendantes.
  
  Soit $v\in V-U$, la courbe $\phi_t^*\cC_v$ dans $\cN_t$ converge quand
  $t\to0$ vers une courbe rationnelle dans $\cN_0$, donc de fibré
  normal $\cO(1)\oplus\cO(1)$.  On en déduit que pour $t\leq t_0$, il en est
  de même pour toutes les courbes $\phi_t^*\cC_v$, donc pour les courbes
  $\cC_{\phi_t(v)}$. Si on prend une famille compacte de vecteurs
  $v$, on peut choisir le même $t_0$ pour tous ces vecteurs.
  
  Si on pouvait appliquer ce résultat pour tous les $v$ au bord d'un
  petit ouvert $B$ de $V$, alors tous les $P^1$ paramétrés par
  $\phi_{t_0}(B)$ auraient fibré normal $\cO(1)\oplus\cO(1)$, et on aurait
  fini.  Malheureusement, la restriction $v\not\in U$ ne permet pas ce
  raisonnement de compacité. Dans un premier temps, on peut déduire
  l'existence de petits disques holomorphes $\Delta$ passant par l'origine
  dans $V$, transverses à $U$, paramétrant des $P^1$ à fibré normal
  $\cO(1)\oplus\cO(1)$.
  
  Dans un second temps, observons que ce résultat vaut pour tous les
  points de $U\subset V$.  Près de l'origine, on décrit ainsi un voisinage
  ouvert de $U$, et le lemme est démontré.
\end{proof}

En poussant plus loin les arguments de cette démonstration, il est
plausible qu'on puisse montrer, sans utiliser le lemme \ref{lem:defC},
l'existence des déformations de la courbe singulière $\cC_x$, avec le
fibré normal attendu : l'idée est, dans le processus
d'approximation ci-dessus, de déformer les courbes rationnelles
provenant du modèle de l'espace de Heisenberg (donc à fibré normal
$\cO(1)\oplus\cO(1)$) en des courbes de $\cN$.  J'ai préféré distinguer les
arguments valables en général (lemme \ref{lem:defC}) des arguments
particuliers à la situation géométrique étudiée.

\section{Le voisinage d'une courbe singulière}
\label{sec:vcs}

Pour réaliser la transformée twistorielle inverse, il sera important
de comprendre de manière précise l'espace des twisteurs $\cN$ au
voisinage d'une courbe singulière.

Comme vu en (\ref{eq:f-ctct}), on peut choisir des coordonnées
$(\sigma,u,v)$ dans un voisinage d'un point $x\in X^\setC$, de sorte que la
forme de contact ait la forme standard :
$$ \eta = d\sigma + \frac{1}{2i}(v du - u dv) . $$
Pour la structure CR standard sur le groupe de Heisenberg, les
$(1,0)$-formes sont données par $\theta^1=du$ et les $(0,1)$-formes par
$\theta^{\bar 1}=dv$.  Par un choix correct de coordonnées et de forme de
contact \cite{CheMos74}, on peut faire coïncider une structure CR
quelconque avec le modèle jusqu'à l'ordre 2 au voisinage du point $x$
:
\begin{equation}
\theta^1-du = O_2(\sigma,u,v), \quad \theta^{\bar 1}-dv = O_2(\sigma,u,v),\label{eq:t1du}
\end{equation}
où $O_2$ signifie des termes d'ordre 2 en les variables indiquées.
La structure réelle s'écrit explicitement comme
\begin{equation}
  \label{eq:tau}
  \tau(\sigma,u,v)=(\bar{\sigma},\bar{v},\bar{u}) .
\end{equation}
Paramétrons la courbe $\cC_y^+=P\Omega'_y$ par $z_+ \to \eta+z_+ \theta^1$ et la
courbe $\cC_y^-$ par $z_- \to \eta + z_- \theta^{\bar{1}}$.
Un voisinage de $\cC_x^+$ dans $\cN^+$ est paramétré par les
coordonnées $(\sigma_+,u_+,z_+)$, où $(\sigma_+,u_+)$ sont les coordonnées
$(\sigma,u)$ sur $\{v=0\}$, et de même un voisinage de $\cC_x^-$ dans
$\cN^-$ est paramétré par les coordonnées $(\sigma_-,v_-,z_-)$.
\begin{lemm}\label{lem:injXN}
  Les injections $X^\setC\to\cN^±$ sont données, dans le cas où $X$ est
  le groupe de Heisenberg, par
  \begin{align*}
    i_+(\sigma,u,v)&=\big(\sigma_+=\sigma+\frac{i}{2}uv,u_+=u,z_+=-iv\big) , \\
    i_-(\sigma,u,v)&=\big(\sigma_-=\sigma-\frac{i}{2}uv,v_-=v,z_-=iu\big) .
  \end{align*}
  Dans le cas général, $i_+$ et $i_-$ diffèrent du modèle par des termes
  d'ordre 3 au moins.
\end{lemm}
\begin{proof}
  Pour l'espace de Heisenberg, on a vu que $\Omega^{1,0}$ est engendré par
  $du$.  Dans les coordonnées $(\sigma,u,v,z_+)$ sur $P\Omega'$, la forme de
  Liouville s'exprime donc comme
\begin{equation}\label{eq:lxpl}
 \lambda = \eta + z_+ du .
\end{equation}
Puisque $d\eta=idu\land dv$, on obtient
\begin{equation}\label{eq:dlxpl}
 d\lambda = (dz_+-idv)\land du .
\end{equation}
Le vecteur $\partial_v+\frac{u}{2i}\partial_\sigma$, de type $(0,1)$, se relève donc en
le vecteur horizontal défini par (\ref{eq:iXdl}),
\begin{equation}\label{eq:XHxpl}
 X_H = \partial_v + \frac{u}{2i}\partial_\sigma  + i\partial_{z_+} .
\end{equation}
Les courbes intégrales de $T^{0,1}$ dans $X^\setC$ sont les courbes
$\{\sigma=\sigma_0,u=u_0\}$, où $\sigma_0$ et $u_0$ sont des constantes ; le vecteur
$X_H$ y étant constant, les feuilles de $P\Omega'$ se décrivent comme
les
$$ F(\sigma_0,u_0,z_0)=\big\{\big(\sigma=\sigma_0 + \frac{uv}{2i},u=u_0,v,z_+=z_0+iv\big)\big\} . $$
En particulier, $(\sigma,u,v,0)$ et $(\sigma-\frac{uv}{2i},u,0,-iv)$ sont dans
la même feuille, ce qui fournit la formule pour $i_+$.

La formule pour $i_-$ est similaire. Pour une variété générale,
l'approximation (\ref{eq:t1du}) indique que $\lambda=\eta+z_+\theta^1$ diffère du
modèle (\ref{eq:lxpl}) à l'ordre 3, donc $d\lambda$ de (\ref{eq:dlxpl}) à
l'ordre 2, ainsi que $X_H$ de (\ref{eq:XHxpl}). Cela donne une
approximation d'ordre 3 sur l'injection $i_+$, obtenue en intégrant $X_H$.
\end{proof}

Les sections de $L$ sur $\cC_y$ s'identifient aux $a+bz_++cz_-$. La
forme de contact tautologique s'exprime alors évidemment dans un
voisinage de $\cC_x$ comme
\begin{align*}
  \eta^c &= d\sigma_+ + z_+ \big(du_+ + O_2(\sigma_+,u_+)\big) ,\\
      &= d\sigma_- + z_- \big(dv_- + O_2(\sigma_-,v_-)\big) .
\end{align*}
En particulier, le long de $\cC_+$, la différentielle $d\eta^c$, bien
définie seulement sur $\ker \eta^c$, est fournie par les formules
$$ d\eta^c|_{\cC_x^+} = dz_+ \land du_+ , \quad
   d\eta^c|_{\cC_x^-} = dz_- \land dv_- . $$
Bien entendu, au point double les deux composantes fournissent le même
résultat $idu\land dv$.

Comme on a vu, il n'y a pas d'obstruction aux déformations de la
courbe rationnelle $\cC_x$, donc l'espace tangent à l'espace des
déformations est paramétré par les sections sur $\cC_x$ du faisceau normal
$(\cI_{\cC_x}/\cI_{\cC_x}^2)^*$, voir \cite[I.2, théorème 2.8]{Kol96}.
En termes concrets, pour le modèle de l'espace de Heisenberg, la
courbe $\cC_x$ a pour équation locale $uv=0$ en coordonnées $(\sigma,u,v)$,
et les sections locales de $\cI_{\cC_x}/\cI_{\cC_x}^2$ sont $d\sigma$ et $u
dv + v du$ ; cela signifie qu'une section locale du faisceau normal
$N_{\cC_x}$ est décrite par un couple $(a\partial_u+b\partial_\sigma,a'\partial_v+b'\partial_\sigma)$ de
sections locales sur $\cC_x^+$ et $\cC_x^-$, avec $b$ et $b'$
holomorphes, $a$ et $a'$ méromorphes avec pôle simple en $x$,
satisfaisant les relations $\Res_xa=\Res_xa'$ et $b=b'$ au point
double $x$. Cette description reste valable dans le cas général par le
lemme \ref{lem:injXN}, et, passant aux systèmes de coordonnées
$(\sigma_+,u_+,z_+)$ et $(\sigma_-,v_-,z_-)$, une section \emph{globale} de $N_{\cC_x}$
est décrite par un couple
\begin{equation}
s=(s_+=a \partial_{u_+}+b \partial_{\sigma_+}, s_-=a' \partial_{v_-}+b' \partial_{\sigma_-} )\label{eq:def-s}
\end{equation}
sur chaque composante, avec $a$ et $a'$ méromorphes sur $\cC_x^+$ et
$\cC_x^-$ avec pôle simple au point double $x$, et $b$ et $b'$
holomorphes (donc constants) ; et satisfaisant la compatibilité au
point double :
\begin{equation}
  \label{eq:compd}
  \Res_x a=-\Res_x a', \quad b'-b=\Res_xa.
\end{equation}
Les sections holomorphes représentent les déformations avec point
double, tandis que les sections méromorphes représentent les lissifications.
Dans la suite, nous utiliserons la base privilégiée de sections
suivantes :
\begin{equation}
  \label{eq:base-s}
  \begin{split}
    s_0&=\frac{1}{i}\big(\frac{1}{z_+}\partial_{u_+}-\frac{1}{2}\partial_{\sigma_+},-\frac{1}{z_-}\partial_{v_-}+\frac{1}{2}\partial_{\sigma_-}\big),\\
    s_1&=\frac{1}{2}(\partial_{\sigma_+},\partial_{\sigma_-}),\\
    s_2&=(\partial_{u_+},0),\\
    s_3&=(0,\partial_{v_-}).
  \end{split}
\end{equation}
On notera $(s^0,s^1,s^2,s^3)$ la base duale. On peut remarquer que la
section $s_2$, tangente à $X^\setC$, engendre $T^{1,0}$, tandis que $s_3$
engendre $T^{0,1}$. En fait, puisqu'on récupère $X^\setC$ comme l'espace
des points doubles des courbes singulières, les vecteurs tangents à
$X^\setC$ en $x$, correspondant respectivement aux sections $s_1$, $s_2$ et
$s_3$, sont en coordonnées $(\sigma,u,v)$ les vecteurs $\frac{\partial_\sigma}{2}$,
$\partial_u$ et $\partial_v$. De (\ref{eq:t1du}) résulte l'égalité
\begin{equation}
  \label{eq:s1s2s3}
  s^1=2\eta, \quad s^2=\theta^1, \quad s^3=\theta^{\bar 1} .
\end{equation}
Enfin, il est facile de voir que la section $s_0$ est réelle (c'est une
raison du facteur $i$ dans son expression).

\section{La transformée twistorielle inverse}\label{sec:la-transf-twist}

Partons d'une variété CR strictement pseudoconvexe $X^3$. On a
construit un espace de twisteurs $\cN$, avec une famille de courbes
rationnelles, à fibré normal $\cO(1)\oplus\cO(1)$, paramétrée par une
variété $M$ de dimension $4$, dans laquelle $X^\setC$ se plonge comme
le diviseur des courbes à points doubles.

L'espace $\cN$ est équipé d'une forme de contact holomorphe
$\eta^c$, à valeurs dans le fibré $L$. Au-dessus d'une courbe singulière
$\cC_x$, le fibré $L$ a degré $1$ sur chaque composante $\cC_x^±$,
donc est au total de degré 2 sur chaque courbe rationnelle.

Dans ces conditions, la transformée twistorielle inverse indique que
$M$ est munie d'une métrique autoduale d'Einstein $g$, dont $\cN$ est
l'espace des twisteurs. Nous emploierons la description suivante de
la métrique.

Considérons le fibré en droites sur $M$ défini par
\begin{equation}
  \label{eq:ell}
  \cL_m = H^0(\cC_m,\Omega^1_{\cC_m}\otimes L) .
\end{equation}
Sur $X\subset M$, la formule continue à avoir un sens, car sur chaque
composante, $L$ se restreint en un faisceau de degré 1, et
$\Omega^1_{\cC_m}$  en un faisceau de degré $-1$.

La forme de contact $\eta^c$, restreinte à une courbe rationnelle
$\cC_x$, définit un élément $\Theta_x\in \cL_x$, donc une section $\Theta$ de
$\cL$. Les courbes à points doubles sont legendriennes pour $\eta^c$,
donc la section $\Theta$ est nulle sur $X$, transverse à $X$ comme on le
verra dans le lemme \ref{lem:gamma}.

Définissons le fibré de rang 3 sur $M$ par $W_m=H^0(\cC_m,L)$, c'est
le fibré $\Omega^2_-M$ des 2-formes antiautoduales sur $M$.  Or l'espace
$H^0(P^1,L)$, où $L=\cO(2)$, est muni d'une métrique canonique $\Upsilon$,
définie par
\begin{equation}
  \label{eq:gammauv}
  \Upsilon(u,v) = u d^2 v + v d^2 u - du dv \in
  H^0(P^1,(\Omega^1_{P^1})^2\otimes L^2)\simeq \setC .
\end{equation}
En formule, si on a une coordonnée $z$ sur $P^1$, et les sections $u$
et $v$ sont données par $u=a+bz+cz^2$ et $v=a'+b'z+c'z^2$, alors
\begin{equation}
  \label{eq:gam-explct}
  \Upsilon(u,v)=\big( 2(ac'+a'c)-bb' \big) dz^2 .
\end{equation}
Cela signifie que le fibré $W$ est muni d'une métrique $\Upsilon$ à valeurs
dans $\cL^2$, donc
$$ \Theta^{-2}\Upsilon $$
est une métrique sur $W$, à savoir (à une constante près) le produit
extérieur sur $\Omega^2_-M$ (défini négatif sur les sections réelles).

On réalise $W_m$ comme espace de 2-formes sur $H^0(\cC_m,N_m)$ de la
manière suivante : si $m\not\in X$, la courbe $\cC_m$ est transverse à
la distribution de contact $\eta^c$, donc le fibré normal $N_m$ de
$\cC_m$ s'identifie à $\ker \eta^c$, et $d\eta^c$ est bien défini sur $N_m$.
On dispose alors d'une 2-forme $\omega$ sur $M-X$, à valeurs dans $W$,
définie par
$$ \omega:\Lambda^2 H^0(\cC_m,N_m) \overset{d\eta^c}{\longrightarrow} W_m . $$
Choisissant une base orthonormale locale $(w_1,w_2,w_3)$ de $W$, nous
obtenons trois 2-formes 
\begin{equation}
\omega_i=2\Theta^{-2}\Upsilon(\omega,w_i)\label{eq:omega-i}
\end{equation}
locales sur $M$ qui forment une base orthonormale de $\Omega^2_-M$. La
connaissance d'une telle base suffit à caractériser la métrique $g$ de
$M$. (Le facteur 2 dans la formule est présent pour normaliser la
métrique de manière agréable dans la suite.)

\begin{rema}
  Nous n'utilisons pas la description, plus habituelle, de la métrique
  comme $g=\omega_H\omega_E$, avec $H_m=H^0(\cC_m,L^{1/2})$,
  $E_m=H^0(\cC_m,N_m\otimes L^{-1/2})$, et $\omega_H$ et $\omega_E$ sont des formes
  symplectiques définies sur $H$ par le wronskien et sur $E$ par
  $d\eta^c$. En effet, la définition de $L^{1/2}$ sur les courbes
  singulières est problématique, car $L$ est de degré 1 sur chaque
  composante. On notera que la métrique (\ref{eq:gammauv}) sur $W$ est
  le carré du wronskien.
\end{rema}

Les sections du faisceau $\Omega^1_{\cC_x}$
s'identifient aux couples $(\alpha_+,\alpha_-)$ de 1-formes sur $\cC_x^+$ et
$\cC_x^-$, méromorphes avec pôle simple au point double $x$, satisfaisant
$$ \Res_x \alpha_+ + \Res_x \alpha_- = 0 . $$
(Cette relation provient de l'équation locale $z_+z_-=0$ pour
$\cC_x$.) Dans la suite, il sera commode de choisir la section
particulière de $\Omega^1_{\cC_x}\otimes L$ sur $\cC_x$, donnée par le couple
$$ \big( \frac{dz_+}{iz_+}\otimes1 , -\frac{dz_-}{iz_-}\otimes1 \big) $$
de sections sur $\cC_x^+$ et sur $\cC_x^-$, pour trivialiser $\cL$ au
point $x$ ; nous étendons cette trivialisation en une trivialisation
$$ \ell : \cL \to \setC $$
dans un voisinage de $x$. En particulier, $\ell\Theta$ devient une
équation locale de $X$ dans $M$, comme le montre le lemme suivant qui
calcule de plus la métrique $\Upsilon$ sur $X$.

\begin{lemm}\label{lem:gamma}
  La section $\Theta$ est transverse sur $X$, donc pour $m\not\in X$, la
  courbe $\cC_m$ est transverse à la distribution de contact $\ker
  \eta^c$.
  
  La métrique $\Upsilon$ sur $W$, à valeurs dans $\cL^2$, satisfait
$$  \ell^2\Upsilon(a+bz_++cz_-,a'+b'z_++c'z_-)=-aa'-2\ell \Theta(bc'+b'c)+O(\Theta^2) .$$
\end{lemm}
\begin{proof}
  Prenons des coordonnées $(z_+,z_-)$ sur $\setC^2$, et définissons
la famille $C_\epsilon$ de courbes de degré 2 de $P^2$ par l'équation
$$ z_+z_-=\epsilon . $$
Pour $\epsilon=0$ on a une courbe singulière $C_0$ et un morphisme évident
(en coordonnées) $f_0:C_0 \to \cC_x \subset \cN$. Celui-ci s'étend en une
famille de morphismes $f_\epsilon:C_\epsilon \to \cN$ induisant une application $\phi$
du disque $\Delta$ dans $M$ telle que $\phi(0)=x$. De plus, on peut
s'assurer que $\frac{d\phi}{d\epsilon}(0)$ est exactement la direction dans $T_xM$
donnée par la section normale $\sigma_0$ définie en (\ref{eq:base-s}).
Alors, sur $\cC_x^+$,
$$ \left.\frac{d\ell \Theta}{d\epsilon}\right|_{\epsilon=0}
 =i_{\sigma_0}d\eta^c(iz_+\partial_{z_+})=d\eta^c(\sigma_0,iz_+ \partial_{z_+})=-1, $$
et on trouve, comme il se doit, le même résultat sur $\cC_x^-$ :
$$ \left.\frac{d\ell \Theta}{d\epsilon}\right|_{\epsilon=0}
 =i_{\sigma_0}d\eta^c(-iz_-\partial_{z_-})=d\eta^c(\sigma_0,-iz_- \partial_{z_-})=-1, $$
d'où 
\begin{equation}
\ell\Theta=-\epsilon + O(\epsilon^2).\label{eq:lT}
\end{equation}
Cela assure bien la transversalité de $\Theta$ sur $X$.

La section $a+bz_++cz_-$ de $L$ induit une section $a+bz_++cz_-$ de
$\cO_{P^2}(1)$, et donc sur chaque $C_\epsilon$ une section
$a+bz_++\frac{c\epsilon}{z_+}$ de $\cO(2)$ (convergeant vers $a+bz_+$ sur
$\cC_x^+$ quand $\epsilon\to0$). En appliquant la formule
(\ref{eq:gam-explct}), on déduit
$$ \ell^2 \Upsilon(a+bz_++cz_-,a'+b'z_++c'z_-)=-aa'+2\epsilon(bc'+b'c)+O(\epsilon^2), $$
et le lemme découle de (\ref{eq:lT}).
\end{proof}

Comme $\cC_x$ est legendrienne pour $\eta^c$, la forme $d\eta^c$ n'est pas
intrinsèque sur le fibré normal de $\cC_x$, mais, pour $s,t\in
H^0(\cC_x,N_x)$ est bien défini
$$ \Theta d\eta^c(s,t)=\eta^c\land d\eta^c(s,t) \in H^0(\cC_x,\Omega^1_{\cC_x}\otimes L^2)=\cL_x\otimes W_x. $$
En calculant cette expression sur la base (\ref{eq:base-s}), on
obtient :
$$ 2\ell\Theta d\eta^c = s^0\land s^1 + z_+ (s^0-is^1)\land s^2 + z_- (s^0+is^1)\land s^3 . $$

Par le lemme \ref{lem:gamma}, une base orthonormale (négative) de $W_x$ pour la
métrique $\Theta^{-2}\Upsilon$ près de $x$, à l'ordre 1 près transversalement à $X$,
est donnée par
\begin{align*}
  w_1 &= \ell \Theta, \\
  w_2 &= \frac{z_++z_-}{2}\sqrt{\ell \Theta}, \\
  w_3 &= \frac{z_+-z_-}{2i}\sqrt{\ell \Theta}.
\end{align*}
Ici il est important de noter que, puisque la structure réelle agit
par $\tau(\eta+z_+\theta^1)=\eta+\bar{z}_+\theta^{\bar 1}$, elle envoie $z_+$ sur
$\bar{z}_-$, si bien que la base ci-dessus est réelle. On en déduit,
par la formule (\ref{eq:omega-i}), la base de 2-formes antiautoduales,
toujours à l'ordre 1 près,
\begin{align*}
  \omega_1 &= (\ell \Theta)^{-2} s^0\land s^1, \\
  \omega_2 &= (\ell \Theta)^{-3/2} \big(
  -(s^0+is^1)\land s^3 -(s^0-is^1)\land s^2 \big),\\
  \omega_3 &= (\ell \Theta)^{-3/2} i \big(
  (s^0+is^1)\land s^3 - (s^0-is^1)\land s^2 \big), \\
\intertext{puis, en posant $\alpha^2=s^2+s^3$ et $\alpha^3=i(s^3-s^2)$, toujours à
$O(\Theta)$ près,}
  \omega_1 &= (\ell \Theta)^{-2} s^0\land s^1, \\
  \omega_2 &= (\ell \Theta)^{-3/2} (-s^0\land\alpha^2-s^1\land\alpha^3), \\
  \omega_3 &= (\ell \Theta)^{-3/2} (s^0\land\alpha^3-s^1\land\alpha^2).  
\end{align*}
Cela correspond aux termes principaux des 2-formes antiautoduales pour
la métrique
$$ \frac{(s^0)^2+(s^1)^2}{(\ell\Theta)^2} + \frac{(\alpha^2)^2+(\alpha^3)^2}{\ell\Theta}
. $$
En revenant à (\ref{eq:s1s2s3}), on a $(\alpha^2)^2+(\alpha^3)^2=4\theta^1\theta^{\bar
  1}=2\gamma$, d'où
$$ g \sim \frac{d(\ell\Theta)^2+4\eta^2}{(\ell\Theta)^2}+\frac{2\gamma}{\ell\Theta}.
$$
Ainsi, la métrique $g$ a près de $X^\setC$ le comportement (\ref{eq:ach}), ce
qui achève la démonstration du théorème \ref{th:1} : suivant qu'on
regarde sur la partie réelle $M^\setR$ le côté $\ell\Theta>0$ ou $\ell\Theta<0$, on
obtient une métrique positive ou de signature $(2,2)$.

Nous serons plus rapides sur la partie unicité du théorème \ref{th:1}.
On a vu en (\ref{eq:courbe-singuliere-x}) que le bord de l'espace des
twisteurs d'une métrique autoduale d'Einstein avec comportement
asymptotique (\ref{eq:ach}) doit être l'espace $\cN|_X$ construit
section \ref{sec:les-twisteurs}, avec sa structure de contact
holomorphe. Le germe de l'espace des twisteurs est uniquement
déterminé par cette hypersurface réelle.

%\bibliographystyle{smfalpha}
%\bibliography{biblio,biquard}

\begin{thebibliography}{AHS78}

\bibitem[AHS78]{AtiHitSin78}
{\scshape M.~F. Atiyah, N.~J. Hitchin {\normalfont \smfandname} I.~M. Singer}
  -- {\og Self-duality in four-dimensional {R}iemannian geometry\fg},
  \emph{Proc. Roy. Soc. London Ser. A} \textbf{362} (1978), no.~1711,
  p.~425--461.

\bibitem[Biq00]{Biq00}
{\scshape O.~Biquard} -- {\og M\'etriques d'{E}instein asymptotiquement
  sym\'etriques\fg}, \emph{Ast\'erisque} \textbf{265} (2000), p.~vi+109.

\bibitem[Biq02]{Biq02}
\bysame , {\og M\'etriques autoduales sur la boule\fg}, \emph{Invent. math.}
  \textbf{148} (2002), no.~3, p.~545--607.

\bibitem[Biq05]{Biq05}
\bysame , {\og Autodual {E}instein versus {K}\"ahler-{E}instein\fg},
  \emph{Geom. Funct. Anal.} \textbf{15} (2005), no.~3, p.~598--633.

\bibitem[Biq06]{Biq06}
\bysame , {\og {Cauchy-Riemann 3-Manifolds and Einstein Fillings}\fg},
  \emph{Perspectives in {R}iemannian {G}eometry} (V.~Apostolov, A.~Dancer,
  N.~Hitchin {\normalfont \smfandname} M.~Wang, \smfedsname), CRM Proceedings
  and Lecture Notes, vol.~40, American Mathematical Society, 2006, p.~27--46.

\bibitem[CM74]{CheMos74}
{\scshape S.~S. Chern {\normalfont \smfandname} J.~K. Moser} -- {\og Real
  hypersurfaces in complex manifolds\fg}, \emph{Acta Math.} \textbf{133}
  (1974), p.~219--271.

\bibitem[CS04]{CalSin04}
{\scshape D.~M.~J. Calderbank {\normalfont \smfandname} M.~A. Singer} -- {\og
  Einstein metrics and complex singularities\fg}, \emph{Invent. Math.}
  \textbf{156} (2004), no.~2, p.~405--443.

\bibitem[Deb01]{Deb01}
{\scshape O.~Debarre} -- \emph{Higher-dimensional algebraic geometry},
  Universitext, Springer-Verlag, New York, 2001.

\bibitem[Fef76]{Fef76}
{\scshape C.~L. Fefferman} -- {\og Monge-{A}mp\`ere equations, the {B}ergman
  kernel, and geometry of pseudoconvex domains\fg}, \emph{Ann. of Math. (2)}
  \textbf{103} (1976), no.~2, p.~395--416.

\bibitem[Fei01]{Fei01}
{\scshape B.~Feix} -- {\og Hyperk\"ahler metrics on cotangent bundles\fg},
  \emph{J. Reine Angew. Math.} \textbf{532} (2001), p.~33--46.

\bibitem[Hit95]{Hit95}
{\scshape N.~J. Hitchin} -- {\og Twistor spaces, {E}instein metrics and
  isomonodromic deformations\fg}, \emph{J. Differential Geom.} \textbf{42}
  (1995), no.~1, p.~30--112.

\bibitem[Kol96]{Kol96}
{\scshape J.~Kollár} -- \emph{Rational curves on algebraic varieties},
  Ergebnisse der Mathematik und ihrer Grenzgebiete., vol.~32, Springer-Verlag,
  Berlin, 1996.

\bibitem[LeB82]{LeB82}
{\scshape C.~LeBrun} -- {\og $\mathcal{H}$-space with a cosmological
  constant\fg}, \emph{Proc. Roy. Soc. London Ser. A} \textbf{380} (1982),
  no.~1778, p.~171--185.

\bibitem[LeB84]{LeB84}
\bysame , {\og Twistor {C}{R} manifolds and three-dimensional conformal
  geometry\fg}, \emph{Trans. Amer. Math. Soc.} \textbf{284} (1984), no.~2,
  p.~601--616.

\bibitem[LM82]{LeeMel82}
{\scshape J.~M. Lee {\normalfont \smfandname} R.~Melrose} -- {\og Boundary
  behaviour of the complex {M}onge-{A}mp\`ere equation\fg}, \emph{Acta Math.}
  \textbf{148} (1982), p.~159--192.

\end{thebibliography}

\providecommand{\bysame}{\leavevmode ---\ }
\providecommand{\og}{``}
\providecommand{\fg}{''}
\providecommand{\smfandname}{et}
\providecommand{\smfedsname}{\'eds.}
\providecommand{\smfedname}{\'ed.}
\providecommand{\smfmastersthesisname}{M\'emoire}
\providecommand{\smfphdthesisname}{Th\`ese}

\end{document}